\documentclass{article}
\usepackage{amsmath}
\usepackage{amsfonts}
\usepackage{amsthm}
\newcommand{\A}{\mathbb{A}}
\newcommand{\RR}{\mathbb{R}}
\newcommand{\ZZ}{\mathbb{Z}}
\newtheorem{theorem}{Theorem}
\newtheorem{lemma}{Lemma}
\newtheorem{cor}{Corollary}
\newtheorem{prop}{Proposition}
\theoremstyle{definition}

\theoremstyle{remark}
\newtheorem*{acn}{Acnowledgements}
\newtheorem*{cau}{Caution}
\DeclareMathOperator{\ord}{ord}
\DeclareMathOperator{\Gal}{Gal}
\DeclareMathOperator{\res}{Res_0}
\DeclareMathOperator{\inj}{in}
\DeclareMathOperator{\md}{md}
\begin{document}
\title{Fractional residues}

\author{Aleksandrs Mihailovs\\
Department of Mathematics\\
University of Pennsylvania\\
Philadelphia, PA 19104-6395\\
mihailov@math.upenn.edu\\
http://www.math.upenn.edu/$\sim$mihailov/
}
\date{\today}
\maketitle

\begin{abstract}
Invariants of generalized tensor fields on a line are classified 
using special polynomials $P_{mk}^{(-1/\lambda)}$ introduced here 
for this purpose. For the case of positive characteristic, 
a new invariant of formal power series, {\em a width}, is defined. 
Some applications to the geometric quantization of a line and 
conformal quantum field theory are discussed as well.
\end{abstract}    
\setlength{\baselineskip}{1.5\baselineskip}

\section{Introduction}\label{sec0}

Differential forms $\phi(t)dt$ on a line have a well-known invariant, 
\begin{equation}\label{in1}
\res \phi(t)dt=a_{-1},
\end{equation}
where
\begin{equation}
\phi(t)=\sum_{i=\ord\phi(t)}^{\infty}a_it^i.
\end{equation}

For quadratic differential forms $\psi(t)(dt)^2$, one can construct an invariant 
by a composition of an invariant mapping 
\begin{equation}
\psi(t)(dt)^2\mapsto \sqrt{\psi(t)} dt
\end{equation} 
and a residue \eqref{in1}. 

Why do we need the invariants of quadratic differential forms? One of the reasons 
is that the pairing 
\begin{equation}
(\psi(t)(dt)^2, \alpha(t)\frac{d}{dt})=\res \psi(t)\alpha(t)dt
\end{equation}
identifies the space of the quadratic differentials with a dual space to the Lie algebra of 
vector fields on a line. Kirillov's orbit method \cite{K,KO1,KO2} associates the orbits 
of the group of authomorphisms of a line in that space with irreducible unitary 
representations of this group. Thus, for geometric quantization of 
a line, we need to describe the orbits and the invariants of quadratic differentials. 

The first calculations for that case were done by the founder of the orbit method, 
Alexandre Kirillov, in \cite{KV}. I was his student at that time, and I found 
the explicit formulas for the invariants, some of them were announced in \cite{KV} 
with an indication of my priority. These results are presented in 
section \ref{sec1}. Most of them are based on the studying of special polynomials 
$P_{mk}$ parametrizing the orbits of the co-adjoint representation of the group of 
authomorphisms of a line.

More generally, the composition of the invariant mapping 
\begin{equation}
\psi(t)(dt)^{-\lambda}\mapsto (\psi(t))^{-1/\lambda} dt
\end{equation}
and a residue \eqref{in1} defines an invariant of generalized differential forms 
$\psi(t)(dt)^{-\lambda}$. This nontrivial invariant allows us to describe the 
orbits of the group of authomorphisms in the space of generalized differential 
forms, utilizing special polynomials $P_{mk}^{(-1/\lambda)}$ which for 
$\lambda=-2$ coinside with polynomials $P_{mk}$ introduced in section \ref{sec1}.
The orbits and invariants of generalized differential forms are described in 
section \ref{sec2}.

Sections \ref{sec1} and \ref{sec2} deal with an arbitrary field of characteristic 0. 
Almost without changes, the results can be transfered to the restricted 
case of a positive characteristic $p>0$. The results related to the geometric quantization 
of a line in the restricted case for $p>0$, corresponding to section \ref{sec1}, 
are described briefly in section \ref{sec3}. 

For the general, not restricted, case of a field $f$ of a positive characteristic 
$p>0$, the situation is much more complicated. In that case even functions have 
a lot of additional invariants. If $f$ is not a perfect field, there are formal 
power series which don't have a polynomial normal form, see \eqref{con}. However, 
all the orbits are closed, and the space of orbits can be metrized by a complete 
metrics. Section \ref{sec4} presents these results. Also, at the end of section 
\ref{sec4} I define a new invariant \eqref{22} of formal vector fields for $p>0$, 
I called it a {\em width}. 

Special polynomials $P_{mk}^{(-1/\lambda)}$ describing the orbits and invariants 
of formal tensor fields on a line, naturally appear in some other fields of 
mathematics as well. An application of them to a particular problem from a quantum 
field theory is discussed in section \ref{sec5}.   

\section{Polynomials $P_{mk}$}\label{sec1}

Let $f$ be an arbitrary field of characteristic 0. Denote $W_1$ the Lie algebra 
of $f$-derivations of $f[[t]]$, the (associative) algebra of the formal power series 
in one variable. In other words, elements of $W_1$ are formal vector fields on 
a line $\A^1(f)=ft$, i.\ e.\ expressions 
\begin{equation}\label{1}
a=\sum_{i=-1}^{\infty}a_il_i,
\end{equation}
with $a_i\in f$ and $l_i=t^{i+1}\frac{\partial}{\partial t}$, with generators 
$l_i$ satisfying 
\begin{equation}\label{2}
[l_i,l_j]=(j-i)l_{i+j} .
\end{equation}

$W_1$ has a natural decreasing filtration 
\begin{equation}\label{3}
W_1=L_{-1}\supset L_0\supset L_1\supset L_2\supset \dots
\end{equation}
where $L_n$ is the Lie subalgebra of $W_1$, consisting of elements \eqref{1} 
with $a_i=0$ if $i<n$. Since $[L_m,L_n]\subseteq L_{m+n}$, Lie algebra $L_n$ is 
an ideal of $L_m$ for all $m$ such that $0\leq m<n$; and we can define Lie 
algebras $L_{mn}=L_m/L_n$ for $0\leq m<n$.

It follows directly from the definition, that $L_{mn}$ is a Lie $f$-algebra of 
dimension $n-m$ with a basis $(l_i+L_n)_{m\leq i<n}$ satisfying 
\begin{equation}\label{4}
[l_i+L_n, l_j+L_n]=\begin{cases} (j-i)(l_{i+j}+L_n) &\text{for $i+j<n$,}\\
0 &\text{otherwise.}\end{cases} 
\end{equation}

Below we'll write $l_i$ instead of $l_i+L_n$, where it won't cause an ambiguity, 
meaning that the brackets $[,]$ in $L_{mn}$ satisfy \eqref{2} with $l_{i+j}=0$ 
for $i+j\geq n$.

Filtration \eqref{3} defines a filtration
\begin{equation}\label{5}
L_{mn}\supset L_{m+1,n}\supset \dots\supset L_{nn}=0
\end{equation}
with
\begin{equation}
[L_{in},L_{jn}]=\begin{cases} L_{i+j,n} &\text{for $i\neq j, i+j<n$,}\\
L_{2i+1,n} &\text{for $i=j, 2i<n$,}\\
0 &\text{otherwise.}\end{cases}
\end{equation}
according to \eqref{4}. Thus $L_{mn}$ are solvable Lie algebras,  
nilpotent for $m>0$ and commutative for $2m+1\geq n$.

Denote $U_{mn}$ the universal enveloping algebra of $L_{mn}$. 
\eqref{5} implies
\begin{equation}
U_{mn}\supset U_{m+1,n}\supset \dots\supset U_{nn}=f
\end{equation}
with $U_{in}=f[l_i,l_{i+1}, \dots, l_{n-1}]$ for $(n-1)/2\leq i\leq n-1$ 
since $L_{in}$ is commutative in these cases. 

\begin{theorem}\label{thm1}
Let $m\geq 0$. If $n\leq 2m+2$, $U_{m+1,n+1}$ is commutative. If $n\geq 2m+2$, 
the center of $U_{m+1,n+1}$ is $f[l_n,l_{n-1},\dots,l_{n-m}]$ for odd $n$ and 
$f[l_n,l_{n-1},\dots,l_{n-m}, P_{mk}(l_n,l_{n-1},\dots,l_{n/2})]$ for even $n$, 
where $k=(n/2)-m$ and polynomials $P_{mk}$ can be defined as the coefficients of 
a generating function 
\begin{equation}\label{6}
\sum_{k=1}^{\infty}P_{mk}(x_0,x_1,\dots,x_{m+k})z^k=
\frac{\sqrt{\sum_{k=0}^{\infty}x_kx_0^{k-1}z^k}-\sqrt{\sum_{k=0}^mx_kx_0^{k-1}z^k}}
{x_0^mz^m}
\end{equation}
\end{theorem}

A canonical projection $L_m\rightarrow L_{mn}$ induces a canonical inclusion of 
the spaces of $f$-linear forms, $L_{mn}^\ast\rightarrow L_m^\ast$, and one has  
an infinite flag
\begin{equation}
0=L_{mm}^\ast\subset L_{m,m+1}^\ast\subset L_{m,m+2}^\ast\subset \dots\subset
L_m^\ast
\end{equation}
with $L_m^\ast=\bigcup_k L_{m,m+k}^\ast$.

Denote $G_0=\Gal(f((t))/f)$, the group of the automorphisms of $f((t))$, the 
field of power series, leaving the constants stable. Elements $g\in G_0$ can be 
uniquely determined by their values $g(t)=ct+o(t)$ with $c\in f^\ast$ i.\ e.\ 
$c\neq 0$. 

For a positive integer $n$, denote $G_n$ a subgroup of $G_0$ consisting of 
elements $g$ satisfying $g(t)=t+o(t^n)$. For $m\leq n$, $G_n$ is a normal subgroup 
of $G_m$. Denote $G_{mn}=G_m/G_n$. 

For $m>0$, there is a standard isomorphism between $G_{mn}$ and an 
(algebraic) adjoint group of a nilpotent Lie algebra $L_{mn}$. The standard 
action of $G_{mn}$ on $L_{mn}^\ast$ coincides with the co-adjoint representation. 
The group $G_m$ also may be considered as an adjoint group of a pronilpotent 
Lie algebra $L_m$ for $m>0$, and the standard action of $G_m$ on $L_m^\ast$ may 
be called a co-adjoint representation. 

\begin{theorem}\label{thm2}
Let $m\geq 0$. If $n\leq 2m+2$, all the orbits of a co-adjoint representation of 
$G_{m+1,n+1}$ in $L_{m+1,n+1}$ are points, and every point is an orbit; all orbits 
are in a general position. If $n>2m+2$, the orbits in a general position of a 
co-adjoint representation $G_{m+1,n+1}$ in $L_{m+1,n+1}$, can be parametrized by 
$m+1$ numbers $c_0\in f^\ast,\thinspace c_1,\dots,c_m\in f$ for odd $n$: they are 
affine planes of dimension $n-2m-1$ defined by equations $l_n=c_0, l_{n-1}=c_1, 
\dots, l_{n-m}=c_m$; or by $m+2$ numbers $c_0\in f^\ast,\thinspace c_1,\dots,c_{m+1}\in f$ 
for even $n$, in which case they are affine varieties of dimension $n-2m-2$ defined 
by equations $l_n=c_0, l_{n-1}=c_1, \dots, l_{n-m}=c_m, P_{mk}(l_n, l_{n-1}, \dots, 
l_{n/2})=c_{m+1}$ with $k=(n/2)-m$ and $P_{mk}$ defined by \eqref{6}. Each orbit 
in a general position of a co-adjoint representation of $G_{m+1,n+1}$ in 
$L_{m+1,n+1}^\ast$, is an orbit of a co-adjoint representation of $G_{m+1,i+1}$ 
in $L_{m+1,i+1}^\ast$ for all $i\geq n$, as well as of $G_{m+1}$ in $L_{m+1}^\ast$.
Each orbit of a co-adjoint representation of $G_{m+1,i+1}$ in $L_{m+1,i+1}^\ast$, 
is an orbit in a general position of a co-adjoint representation of $G_{m+1,n+1}$ in 
$L_{m+1,n+1}^\ast$ for some $n$ such that $\min (2m+2,i)\leq n\leq i$. Each orbit of 
a co-adjoint representation of $G_{m+1}$ in $L_{m+1}^\ast$ is an orbit in a general 
position of a co-adjoint representation of $G_{m+1,n+1}$ in $L_{m+1,n+1}^\ast$ 
for some $n\geq 2m+2$.
\end{theorem}

Let us study polynomials $P_{mk}$ mentioned in Theorems \ref{thm1} and \ref{thm2} 
in more details. 

\begin{theorem}\label{thm3}
A polynomial $P_{mk}$ is homogeneous of degree $k$ and equalized of weight $m+k$.
\begin{equation}\label{7}
P_{mk}(x_0,\dots,x_{m+k})=\sum_{\substack{\pi\vdash m+k\\ \pi_1>m}}
\binom{1/2}{p_1,\dots,p_{m+k}}x_0^{p_0}x_1^{p_1}\dots x_{m+k}^{p_{m+k}} ,
\end{equation}
where $\pi=(1^{p_1}2^{p_2}\dots)$ is supposed to be a partition of $m+k$ with the  
largest part $\pi_1>m$; $p_0=k-\ell (\pi)$ where $\ell (\pi)$ denotes the length of 
a partition $\pi$. The least common multiple of the denominators of the coefficients 
of $P_{mk}$ equals $2^{2k-s(k)}$ for $m=0$ or $2^{2k-s(k-1)-1}$ for $m>0$, where 
$s(k)$ is the sum of digits of the binary expression of $k$. Also, 
\begin{equation}\label{8}
P_{mk}(x_0,\dots,x_{m+k})=\frac{(-1)^{k-1}}{(2k-2)!!\cdot 2}\thickspace 
\int_0^x\det\binom{dx}{A} ,
\end{equation}
where $x=(x_{m+1},\dots,x_{m+k}),\thinspace dx=(dx_{m+1},\dots,dx_{m+k})$ and 
\begin{equation}\label{9}
A=\begin{pmatrix}(2k-2)x_0 &(2k-3)x_1 &\dots &kx_{k-2} &(k-1)x_{k-1}\\ 
0 &(2k-4)x_0 &\dots &(k-1)x_{k-3} & (k-2)x_{k-2}\\ 
\hdotsfor[2]{5}\\ 0 &0 &\dots &2x_0 &x_1
\end{pmatrix} 
\end{equation}
is a $(k-1)\times k$ matrix; $\binom{dx}{A}$ denotes $k\times k$ matrix obtained 
from $A$ by adding a first row $dx$. Also, 
\begin{equation}\label{10}
P_{0k}(x_0,\dots,x_k)=\frac{(-1)^{k-1}}{(2k)!!}\det\binom{x'}{A} ,
\end{equation}
where $x'=(x_1, 2x_2,\dots,kx_k)$ and $A$ as in \eqref{9}. If $m\geq k-1$, 
then 
\begin{equation}
P_{mk}(x_0,\dots,x_{m+k})=\frac{(-1)^{k-1}}{(2k-2)!!\cdot 2}\det\binom{x}{A} 
\end{equation}
with $x=(x_{m+1},\dots,x_{m+k})$ and matrix $A$ defined above in \eqref{9}.
One has $P_{m1}=\frac{1}{2}x_{m+1}$;\\ 
$P_{02}=\frac{1}{8}(4x_0x_2-x_1^2)$ and $P_{m2}=\frac{1}{4}(2x_0x_{m+2}-x_1x_{m+1})$ 
for $m\geq 1$;\\
$P_{03}=\frac{1}{16}(8x_0^2x_3-4x_0x_1x_2+x_1^3),\quad 
P_{13}=\frac{1}{16}(8x_0^2x_4-2x_0(2x_1x_3+x_2^2)+3x_1^2x_2)$,\\ 
$P_{m3}=\frac{1}{16}(8x_0^2x_{m+3}-4x_0(x_1x_{m+2}+x_2x_{m+1})+3x_1^2x_{m+1})$ for 
$m\geq 2$;\\ 
$P_{04}=\frac{1}{128}(64x_0^3x_4-16x_0^2(2x_1x_3+x_2^2)+24x_0x_1^2x_2-5x_1^4)$,\\
$P_{14}=\frac{1}{32}(16x_0^3x_5-8x_0^2(x_1x_4+x_2x_3)+6x_0(x_1^2x_3+x_1x_2^2)
-5x_1^3x_2)$,\\ 
$P_{24}=\frac{1}{32}(16x_0^3x_6-4x_0^2(2x_1x_5+2x_2x_4+x_3^2)+6x_0(x_1^2x_4+
2x_1x_2x_3)-5x_1^3x_3)$,\\ 
$P_{m4}=\frac{1}{32}(16x_0^3x_{m+4}-8x_0^2(x_1x_{m+3}+x_2x_{m+2}+x_3x_{m+1})+
6x_0(x_1^2x_{m+2}+2x_1x_2x_{m+1})-5x_1^3x_{m+1})$ for $m\geq 3$. Also, 
\begin{equation}\begin{split}
P_{0k}(1,1,\dots,1)=&\frac{(2k-1)!!}{(2k)!!}  \\
P_{mk}(1,1,\dots,1)=&\frac{(2k-3)!!}{(2k-2)!!\cdot 2}\quad \text{for $m\geq k-1$.}
\end{split}\end{equation}
If $m=0$, the sum \eqref{7} has $p(k)$, the number of partitions of $k$, 
nonzero items. If $m\geq k-1$, the sum 
\eqref{7} has $p(0)+p(1)+\dots+p(k-1)$ nonzero items.
\end{theorem}

Referring to $x_0,x_1,\dots,x_m$ as constants, one obtains from \eqref{8}, 
\begin{equation}\label{11}
dP_{mk}=\frac{\partial P_{mk}}{\partial x_{m+1}} dx_{m+1}+\dots+
\frac{\partial P_{mk}}{\partial x_{m+k}} dx_{m+k}=\det\binom{dx}{A} .
\end{equation}
Expanding the determinant along the first row, we get determinant formulas for 
partial derivatives:

\begin{cor}
For an integer $i$ so that $1\leq i\leq m$, 
\begin{equation}
\frac{\partial P_{mk}}{\partial x_{m+i}}=(-1)^{i+1}\det A_i ,
\end{equation} 
where $A_i$ is a matrix obtained from $A$ by deleting $i$-th column.
\end{cor}

\begin{lemma}\label{lem1}
Let $r$ be an arbitrary commutative ring, $\frac{\partial}{\partial x_1},\dots, 
\frac{\partial}{\partial x_k}$ some derivations of $r$, $P\in r$ and 
\begin{equation}
dP \stackrel{\rm def}{=} \frac{\partial P}{\partial x_1}dx_1+\dots+
\frac{\partial P}{\partial x_k}dx_k=\det \begin{pmatrix}
dx_1&\dots&dx_k\\ 
a_{21}&\dots&a_{2k}\\
\hdotsfor[1.5]{3}\\
a_{k1}&\dots&a_{kk}\end{pmatrix}
\end{equation}
with $a_{ij}\in r$. Then $P$ satisfies the following system of partial 
differential equations:
\begin{equation}\label{12}
\left\{ \begin{array}{l}a_{21}\frac{\partial P}{\partial x_1}+\dots+
a_{2k}\frac{\partial P}{\partial x_k}=0 ,\\
\hdotsfor[2]{1}\\
a_{k1}\frac{\partial P}{\partial x_1}+\dots+
a_{kk}\frac{\partial P}{\partial x_k}=0 . 
\end{array}\right.
\end{equation}
\end{lemma}
\begin{proof}
Expanding the determinant 
\begin{equation}
\det \begin{pmatrix} a_{i1} &\dots &a_{ik}\\ 
a_{21} &\dots &a_{2k}\\ 
\hdotsfor[1.5]{3}\\ 
a_{i1} &\dots &a_{ik}\\
\hdotsfor[1.5]{3}\\
a_{k1} &\dots &a_{kk} \end{pmatrix} =0
\end{equation}
along the first row, we get the corresponding equation of the system \eqref{12}.
\end{proof}

\begin{proof}[Proof of Theorem \ref{thm1}] 
By Gelfand's Lemma \cite{Gel}, for the standard representation of a Lie 
algebra having a basis $(x_0,\dots,x_N)$, on the polynomial algebra 
$f[x_0,\dots,x_N]$, one has 
\begin{equation}\label{13}
T(x_i)P=\sum_{j=0}^N [x_i,x_j] \frac{\partial P}{\partial x_j}
\end{equation}
for all $i$ from $0$ to $N$ and $P\in f[x_0,\dots,x_N]$. 
Thus the algebra of invariants of the given Lie algebra in $f[x_0,\dots,x_N]$, 
can be described as the algebra of solutions of a system of partial 
differential equations 
\begin{equation}\label{sys}  
\left\{ \begin{array}{l}
\sum_{j=0}^N [x_0,x_j] \frac{\partial P}{\partial x_j}=0 ,\\
\hdotsfor[1.5]{1}\\
\sum_{j=0}^N [x_N,x_j] \frac{\partial P}{\partial x_j}=0 .
\end{array} \right.
\end{equation} 

For a Lie algebra $L_{m+1,n+1}$ with odd $n$, the system \eqref{sys} where 
$x_i=l_{n-i}$ and $N=n-m-1$, has an upper triangular matrix with a 
non-zero main diagonal. Thus, by induction, 
invariants don't depend on $l_i$ with $m<i<n-m$. If $n$ is even, the matrix 
of the coefficients of system 
\eqref{sys} for $L_{m+1,n+1}$ has a form
\begin{equation}
\begin{pmatrix}A&B\\ 0&C\end{pmatrix}
\end{equation}
where $A$ is defined in \eqref{9}, and $C$ is an upper triangular square 
matrix with a non-zero main diagonal. By induction, the same as for the case 
of an odd $n$, one can deduct that invariants don't depend on $l_i$ for 
$m<i<n/2$. Further, the invariants satisfy the system of partial differential 
equations with a matrix $A$. Now, utilizing Lemma \ref{lem1} and formulas 
\eqref{8}, \eqref{11}, we get an additional invariant $P_{mk}(l_n,\dots,l_{n/2})$,  
with $k=(n/2)-m$.

It follows from the general theory of invariants of nilpotent Lie algebras  
\cite{MD}, that the algebra of invariants of $L_{m+1,n+1}$ 
discussed above, for even $n$, is $f[l_n, l_{n-1},\dots, l_{n-m}, P]$ with unknown 
polynomial $P$. Notice that our polynomial $P_{mk}(l_n, \dots, l_{n/2})$ with 
$k=(n/2) - m$, is not included in any algebras $f[l_n, l_{n-1}, \dots, l_{n-m}, P]
\subseteq f[l_n,\dots,l_{n/2}]$ such that $P \not\in f[l_n, \dots, l_{n-m}, 
P_{mk}(l_n, \dots, l_{n/2})]$, since 
\begin{equation}
P_{mk}(l_n, \dots, l_{n/2})=\frac{1}{2}l_n^{k-1}l_{n/2}+l_nQ_{mk}(l_n, \dots, 
l_{(n/2)+1}) + c_{mk}l_{n-1}^{k-1}l_{n-m-1}
\end{equation} 
for some polynomial $Q_{mk}$ and nonzero constant $c_{mk}$, meaning 
that $P_{mk}$ is a linear polynomial of $l_{n/2}$ with coprime coefficients.

$l_n,\dots,l_{n-m}$ are central elements of $L_{m+1,n+1}$. Connections between 
invariants and the center of $U_{m+1,n+1}$, 
the universal enveloping algebra, are well known now, and can be found in 
\cite{Gel}. \end{proof}

Another proof of Theorem \ref{1}, based on the studying of generating functions 
\eqref{6}, will be given in the next section. Theorem \ref{thm3} (except the determinant 
formula \eqref{10}) follows from a comparison between these two proofs of 
Theorem \ref{thm1}. The determinant formula can be obtained by differentiation of 
the corresponding generating function \eqref{6} and observing the conditions on 
coefficients; similar calculations can be found in \cite{MM} and \cite{WG}.
Theorem \ref{thm2} follows from Theorem \ref{thm1} and the results of \cite{MD}. 

Some of results of this section were announced in \cite{KV}.  

\section{Fractional residues}\label{sec2}

The same as in the previous section, let $f$ be a field of characteristic $0$.
For $\lambda, \mu \in f$ denote $F_{\lambda\mu} = f[[t]]t^\mu(dt)^{-\lambda}$,   
a linear topological $f$-space with the topology induced from the 
standard topology of $f[[t]]$, the algebra of formal power series, assuming a 
discrete topology of $f$. Lie algebras $L_m$ and groups $G_m$, 
with $m>0$, naturally act on these spaces. The purpose of this section is to study 
the algebras $I_{\lambda\mu}^m$ of polynomial invariants of these actions.

Elements $e_k=t^{k+\mu}(dt)^{-\lambda}$, where $k=0, 1, 2, \dots$, form a topological 
basis of $F_{\lambda,\mu}$. Denote $(x_k)_{0\leq k\in \ZZ}$ the dual basis of the 
topological $f$-space $F_{\lambda,\mu}^\ast$ of linear forms on 
$F_{\lambda,\mu}$.

\begin{theorem}\label{thm4}
Let $m$ be a non-negative integer. If $\mu \neq (m+k+1)\lambda$ for any positive 
integers $k$, then $I_{\lambda\mu}^{m+1}=f[x_0, x_1, \dots, x_m]$. If 
$\lambda=\mu=0$, then $I_{\lambda\mu}^{m+1}=f[x_0, x_1, \dots, x_m+1]$. If 
$\lambda\neq 0, \mu=(m+k+1)\lambda$ for a positive integer $k$, and $-1/\lambda
\neq n$ for any positive integer $n<k$, then $I_{\lambda\mu}^{m+1} = 
f[x_0, x_1, \dots, x_m, P_{mk}^{(-1/\lambda)}(x_0, x_1, \dots, x_{m+k})]$ where 
polynomial $P_{mk}^{(-1/\lambda)}$ is defined by a generating function 
\begin{equation}\label{14}
\sum_{k=1}^{\infty} P_{mk}^{(-1/\lambda)}(x_0, x_1, \dots, x_{m+k}) z^k = 
\frac{(\sum_{i=0}^{\infty}x_ix_0^{i-1}z^i)^{-1/\lambda} - 
(\sum_{i=0}^{m}x_ix_0^{i-1}z^i)^{-1/\lambda}}{x_0^mz^m} .
\end{equation}
If $-1/\lambda=n$ for a positive integer $n$ and $\mu=(m+k+1)\lambda$ for a positive
integer $k>n$, then $I_{\lambda\mu}^{m+1} = f[x_0, x_1, \dots, x_m, 
P_{mk}^{(n)}(x_0, x_1, \dots, x_{m+k})/x_0^{k-n}]$ 
where polynomial $P_{mk}^{(n)}$ is defined above.
\end{theorem}

\begin{theorem}\label{thm5}
Let $m$ be a non-negative integer. If $\mu \neq (m+k+1)\lambda$ for any positive 
integers $k$, then the orbits in a general position of the standard representation 
of $G_{m+1}$ in $F_{\lambda\mu}$ can be parametrized by $(m+1)$ numbers 
$c_0\in f^\ast, c_1, \dots, c_m\in f$: they are affine planes of codimension $m+1$  
given by equations $x_0=c_0, x_1=c_1, \dots, x_m=c_m$. If 
$\lambda\neq 0, \mu=(m+k+1)\lambda$ for a positive integer $k$, then the orbits 
in a general position of the standard representation of  
$G_{m+1}$ in $F_{\lambda\mu}$ can be parametrized by $(m+2)$ numbers 
$c_0\in f^\ast, c_1, \dots, c_m, c_{m+1}\in f$: they are affine varieties of 
codimension $m+2$ given by equations $x_0=c_0, x_1=c_1, \dots, x_m=c_m, 
P_{mk}^{(-1/\lambda)}(x_0, x_1, \dots, x_{m+k})=c_{m+1}$ if $-1/\lambda\neq n$ 
for any positive integer $n<k$, or $x_0=c_0, x_1=c_1, \dots, x_m=c_m, 
P_{mk}^{(-1/\lambda)}(x_0, x_1, \dots, x_{m+k})/x_0^{k-n}=c_{m+1}$ if 
$-1/\lambda=n$ for a positive integer $n$, where polynomial $P_{mk}^{(-1/\lambda)}$ 
is defined by \eqref{14}. Each orbit in a general position of the standard 
representation of $G_{m+1}$ in $F_{\lambda\mu}$ is an orbit of the standard 
representation of $G_{m+1}$ in $F_{\lambda,\mu-i}$ for each nonnegative integer 
$i$. Each orbit of the standard representation of $G_{m+1}$ in $F_{\lambda\mu}$ is
an orbit in a general position of the standard 
representation of $G_{m+1}$ in $F_{\lambda,\mu+i}$ for a nonnegative integer $i$, 
with the only exception when $\lambda=0$ and $\mu$ is a non-positive 
integer: then sets $c+{\cal O}$ are also orbits for any $c\in f^\ast$ and $\cal O$, 
an orbit in general position of the standard 
representation of $G_{m+1}$ in $F_{0i}$ for a positive integer $i$.
\end{theorem}

Let us study polynomials $P_{mk}^{(-1/\lambda)}$ mentioned in Theorems \ref{thm4} 
and \ref{thm5} in more details. 

\begin{theorem}\label{thm6}
A polynomial $P_{mk}^{(-1/\lambda)}$ 
is homogeneous of degree $k$ and equalized of weight $m+k$.
\begin{equation}\label{15}
P_{mk}^{(-1/\lambda)}(x_0,\dots,x_{m+k})=\sum_{\substack{\pi\vdash m+k\\ \pi_1>m}}
\binom{1/\lambda}{p_1,\dots,p_{m+k}}x_0^{p_0}x_1^{p_1}\dots x_{m+k}^{p_{m+k}} ,
\end{equation}
where $\pi=(1^{p_1}2^{p_2}\dots)$ is supposed to be a partition of $m+k$ with the  
largest part $\pi_1>m$; $p_0=k-\ell (\pi)$ where $\ell (\pi)$ denotes the length of 
a partition $\pi$. Also, 
\begin{equation}
P_{mk}^{(-1/\lambda)}(x_0,\dots,x_{m+k})=\frac{1}{(k-1)! (-\lambda)^k}\thickspace 
\int_0^x\det\binom{dx}{A} ,
\end{equation}
where $x=(x_{m+1},\dots,x_{m+k}),\thinspace dx=(dx_{m+1},\dots,dx_{m+k})$ and 
\begin{equation}\label{16}
A=\begin{pmatrix}(k-1)\lambda x_0 &((k-1)\lambda+1)x_1 &\dots &((k-1)\lambda+(k-2))x_{k-2} &(k-1)(\lambda+1)x_{k-1}\\ 
0 &(k-2)\lambda x_0 &\dots &((k-2)\lambda+(k-3))x_{k-3} & (k-2)(\lambda+1)x_{k-2}\\ 
\hdotsfor[2]{5}\\ 0 &0 &\dots &\lambda x_0 &(\lambda+1)x_1
\end{pmatrix} 
\end{equation}
is a $(k-1)\times k$ matrix; $\binom{dx}{A}$ denotes $k\times k$ matrix obtained 
from $A$ by adding a first row $dx$. Also, 
\begin{equation}
P_{0k}^{(-1/\lambda)}(x_0,\dots,x_k)=\frac{1}{k! (-\lambda)^k}\det\binom{x'}{A} ,
\end{equation}
where $x'=(x_1, 2x_2,\dots,kx_k)$ and $A$ as in \eqref{16}. If $m\geq k-1$, 
then 
\begin{equation}
P_{mk}^{(-1/\lambda)}(x_0,\dots,x_{m+k})=\frac{1}{(k-1)! (-\lambda)^k}\det\binom{x}{A} 
\end{equation}
with $x=(x_{m+1},\dots,x_{m+k})$ and matrix $A$ defined above in \eqref{16}.
One has $P_{m1}^{(-1/\lambda)}=-\frac{1}{\lambda}x_{m+1}$;\\ 
$P_{02}^{(-1/\lambda)}=-\frac{1}{\lambda}x_0x_2+\frac{\lambda+1}{2\lambda^2}x_1^2$ 
and $P_{m2}{(-1/\lambda)}=-\frac{1}{\lambda}x_0x_{m+2}+\frac{\lambda+1}{\lambda^2}x_1x_{m+1}$ 
for $m\geq 1$;\\
$P_{03}^{(-1/\lambda)}=-\frac{1}{\lambda}x_0^2x_3+\frac{\lambda+1}{\lambda^2}x_0x_1x_2-
\frac{(\lambda+1)(2\lambda+1)}{6\lambda^3}x_1^3,\\
P_{13}^{(-1/\lambda)}=-\frac{1}{\lambda}x_0^2x_4+\frac{\lambda+1}{2\lambda^2}x_0(2x_1x_3+x_2^2)-
\frac{(\lambda+1)(2\lambda+1)}{2\lambda^3}x_1^2x_2$,\\ 
$P_{m3}^{(-1/\lambda)}=-\frac{1}{\lambda}x_0^2x_{m+3}+\frac{\lambda+1}{\lambda^2}x_0(x_1x_{m+2}+x_2x_{m+1})- 
\frac{(\lambda+1)(2\lambda+1)}{2\lambda^3}x_1^2x_{m+1}$ for $m\geq 2$;\\ 
$P_{04}^{(-1/\lambda)}=-\frac{1}{\lambda}x_0^3x_4+\frac{\lambda+1}{2\lambda^2}x_0^2(2x_1x_3+x_2^2)-
\frac{(\lambda+1)(2\lambda+1)}{2\lambda^3}x_0x_1^2x_2+
\frac{(\lambda+1)(2\lambda+1)(3\lambda+1)}{24\lambda^4}x_1^4$,\\
$P_{14}^{(-1/\lambda)}=-\frac{1}{\lambda}x_0^3x_5+\frac{\lambda+1}{\lambda^2}x_0^2(x_1x_4+x_2x_3)-
\frac{(\lambda+1)(2\lambda+1)}{2\lambda^3}x_0(x_1^2x_3+x_1x_2^2)+
\frac{(\lambda+1)(2\lambda+1)(3\lambda+1)}{6\lambda^4}x_1^3x_2$,\\ 
$P_{24}^{(-1/\lambda)}=-\frac{1}{\lambda}x_0^3x_6+\frac{\lambda+1}{2\lambda^2}x_0^2(2x_1x_5+2x_2x_4+x_3^2)-
\frac{(\lambda+1)(2\lambda+1)}{2\lambda^3}x_0(x_1^2x_4+
2x_1x_2x_3)+\frac{(\lambda+1)(2\lambda+1)(3\lambda+1)}{6\lambda^4}x_1^3x_3$,\\ 
$P_{m4}^{(-1/\lambda)}=-\frac{1}{\lambda}x_0^3x_{m+4}+
\frac{\lambda+1}{\lambda^2}x_0^2(x_1x_{m+3}+x_2x_{m+2}+x_3x_{m+1})-
\frac{(\lambda+1)(2\lambda+1)}{2\lambda^3}x_0(x_1^2x_{m+2}+2x_1x_2x_{m+1})+
\frac{(\lambda+1)(2\lambda+1)(3\lambda+1)}{6\lambda^4}x_1^3x_{m+1}$ for $m\geq 3$. Also, 
\begin{equation}\begin{split}
P_{0k}^{(-1/\lambda)}(1,1,\dots,1)=&(-1)^k\binom{1/\lambda}{k} ,  \\
P_{mk}^{(-1/\lambda)}(1,1,\dots,1)=&\frac{(-1)^k}{\lambda}
\binom{1/\lambda}{k-1}\quad \text{for $m\geq k-1$.}
\end{split}\end{equation}
If $-1/\lambda \neq n$ for any positive integer $n<k$,  
the sum \eqref{15} has $p(k)$  
nonzero items for $m=0$, or $p(0)+p(1)+\dots+p(k-1)$ nonzero items 
for $m\geq k-1$. If $-1/\lambda = n$ for a positive integer $n<k$, then 
the sum \eqref{15} has $p_n(k)$, the number of partitions of $k$ with 
length $\leq n$, nonzero items for $m=0$, or $p_n(0)+p_n(1)+\dots+p_n(k-1)$ 
nonzero items for $m\geq k-1$. 
\end{theorem}

\begin{proof}[Proof of Theorem \ref{thm4}] 
Analogously to the proof of Theorem \ref{1}, one can notice that for a representation $T$ 
of a Lie $f$-algebra with a basis $(l_i)_{i\in J}$, in the algebra 
$f[x_0, x_1, \dots]$, one has the series of equalities 
\begin{equation}
T(l_i)P=\sum_{j\in J}(l_ix_j)\frac{\partial P}{\partial x_j} .
\end{equation}
Continuing as in the proof of Theorem \ref{thm1}, we get a proof of Theorem 
\ref{thm4}.

There is also another proof. Consider a residue $\res=x_k\in I_{-1,-k-1}$ where 
$k$ is a nonnegative integer. For such a non-negative $k$, if $\mu=(k+1)\lambda$ 
and $\lambda \neq 0$, then one has an invariant polynomial mapping 
\begin{equation}\label{17}\begin{split}
P^{(-1/\lambda)}: \left.F_{\lambda\mu}\right|_{x_0=1}&\longrightarrow 
F_{-1,-k-1}\\ ht^\mu(dt)^{-\lambda}&\mapsto (ht^\mu(dt)^{-\lambda})^{-1/\lambda}=
h^{-1/\lambda}t^{-k-1}dt .   
\end{split}
\end{equation}
Composing this invariant mapping with a standard residue $\res$, we obtain 
a polynomial invariant for every positive integer $k$,
\begin{equation}\label{18}
\res \circ P^{(-1/\lambda)}: \left. F_{\lambda\mu}\right|_{x_0=1} \longrightarrow f .
\end{equation}

Noticing that $x_0$ is also a $G_1$-invariant, we can extend \eqref{18} first to a 
rational $G_1$-invariant 
\begin{equation}\label{19}
\res \circ P^{(-1/\lambda)} \circ\left( \frac{\cdot}{x_0}\right) : 
\left. F_{\lambda\mu}\right|_{x_0\neq 0} \longrightarrow f 
\end{equation}
and then to a polynomial $G_1$-invariant 
\begin{equation}\label{20}
P_{0k}^{(-1/\lambda)}=\left( \cdot x_0^k\right) \circ\res\circ 
P^{(-1/\lambda)} \circ\left( \frac{\cdot}{x_0}\right) :
F_{\lambda\mu}\longrightarrow f .
\end{equation}
The rest of the proof can be done by utilizing the standard techniques from 
\cite{MD}.
\end{proof} 
 
Noticing that bilinear transformations 
\begin{equation}
P: F_{\lambda\mu}\times F_{\lambda'\mu'}\longrightarrow F_{\lambda+\lambda', 
\mu+\mu'}, (e_i, e_j)\mapsto e_{i+j}
\end{equation}
are invariant, one obtains a bilinear invariant for $\lambda+\lambda'=-1$ when 
$\mu+\mu'$ is a negative integer:
\begin{equation}
\res\circ P: F_{\lambda\mu}\times F_{\lambda'\mu'}\longrightarrow f.
\end{equation}
Thus, 
\begin{equation}
F_{\lambda\mu}^\ast\simeq \left( \bigcup_i F_{-1-\lambda,i-\mu}\right) / 
F_{-1-\lambda,-\mu}\stackrel{\rm def}{=}F_{-1-\lambda,-\mu}^- .
\end{equation}

In particular, 
\begin{equation}
L_m^\ast=F_{1,m+1}\simeq F_{-2,-1-m}^-
\end{equation}
and
\begin{equation}
L_{mn}^\ast\simeq F_{-2,-1-n}/F_{-2,-1-m} .
\end{equation}
That explains the identity 
\begin{equation}
P_{mk}=P_{mk}^{(1/2)}
\end{equation}
following from \eqref{6} and \eqref{14}.

\begin{cau}
`An orbit in a general position', here and in the previous chapter, 
doesn't mean the `orbit of maximal dimension' or the 
`orbit of minimal codimension'.  For instance, the following two 
series of the orbits of the co-adjoint representation of $G_{1,5}$: given 
by equations 
\begin{equation}\label{gen}
l_4=c_0,\quad \frac{1}{2}l_4l_2-\frac{1}{8}l_3^2=c_1
\end{equation}
with $c_0\in f^\ast, c_1\in f$, and by equations 
\begin{equation}
l_4=0,\quad l_3=c_0
\end{equation}
with $c_0\in f^\ast$, both have the maximal dimension 2, but only \eqref{gen} 
are orbits in a general position of the co-adjoint representation of $G_{1,5}$
in $L_{1,5}^\ast$.
\end{cau}

\section{Positive characteristic, a restricted case}\label{sec3}

Let $f$ be an arbitrary field of a positive characteristic $p$. 
Denote $W_1$ the restricted Lie $p$-algebra of $f$-derivations of 
$f[t]/(t^p)$. Elements of $W_1$ can be written in a form 
\begin{equation}\label{21}
a=\sum_{i=-1}^{p-2} a_i l_i
\end{equation}
with $a_i\in f$ and $l_i=t^{i+1}\frac{d}{dt}$. The basic elements $l_i$ 
satisfy \eqref{2} meaning $l_{i+j}=0$ for $i+j\geq p-1$. Also 
$l_i^p=0$ for $i\neq 0$, and $l_0^p=l_0$. 

The same as in section \ref{sec1}, consider filtration 
\begin{equation}
W_1=L_{-1}\supset L_0\supset \dots \supset L_{p-2}
\end{equation}
assuming that $L_n$ is a Lie $p$-algebra consisting of expressions \eqref{21} 
with $a_i<0$ for $i<n$. Define $L_{mn}=L_m/L_n$ for $0\leq m\leq n\leq p-1$, 
supposing that $L_{m,p-1}=L_m$. Denote $U_{mn}$ the restricted universal 
enveloping algebra of $L_{mn}$.

\begin{theorem}\label{thm7}
Let $m\geq 0$. If $n\leq 2m+2$, then $U_{m+1,n+1}$ is commutative. If 
$p-2\geq n\geq 2m+2$, the center of $U_{m+1,n+1}$ is 
$f[l_n, l_{n-1}, \dots, l_{n-m}]/(l_n^p, l_{n-1}^p, \dots, l_{n-m}^p)$ for odd $n$, 
or $f[l_n, l_{n-1}, \dots, l_{n-m},\\ P_{mk}(l_n, l_{n-1},\dots, l_{n/2})]/ 
((l_n^p,\dots,l_{n/2}^p)\bigcap 
f[l_n, l_{n-1},\dots, l_{n-m}, P_{mk}(l_n, l_{n-1},\dots, l_{n/2})])$ for even $n$, 
where $k=(n/2)-m$, and $P_{mk}$ defined by \eqref{6}.
\end{theorem}

Denote $G_0$ the group of automorphisms of $f[t]/(t^p)$, leaving the constants 
stable, and $G_n$ with $1\leq n\leq p-1$, the subgroup of $G_0$ of 
automprphisms $g(t)=t+o(t^n)$. Also denote $G_{mn}=G_m/G_n$ for $0\leq m\leq n
\leq p-1$.

\begin{theorem}\label{thm8}
Theorem \ref{thm2} is true mutatis mutandis.
\end{theorem}

The proofs of Theorems \ref{thm7} and \ref{thm8} can be obtained the same way 
as the proofs of Theorems \ref{thm1} and \ref{thm2}, mutatis mutandis. 

\section{Formal singularities}\label{sec4}

Let $f$ be an arbitrary field of characteristic $p\geq 0$. Denote 
$F=f((t))$ and $G_0=\Gal(F/f)$. We suppose that $F$ has a valuation, 
a filtration and a topology, as usual, assuming a discrete topology of $f$. 
One can check that all of the elements of $G_0$ are automatically 
continous automorphisms preserving the valuation, and are defined uniquely 
by their values $g(t)=ct+o(t)$ with $c\in f^\ast$, i.\ e. $c\neq 0$.

For a positive integer $n$, denote $G_n$ the subgroup of $G_0$ of automorphisms 
$g$ such that $g(t)=t+o(t^n)$. The filtration 
\begin{equation} 
G_0\supset G_1\supset G_2\supset \dots
\end{equation} 
defines a topology (`given by a filtration') in every $G_n$. 

\begin{theorem}\label{thm9}
Let $n$ be a nonnegative integer. $G_n$-orbit of a formal power series 
$h\in F$, is open iff $h\in F\setminus f((t^p))$. The statement, 
`$G_n$-orbit of $h\in f((t^{p^m}))\setminus f((t^{p^{m+1}}))$ is open in the 
relative topology of a field $f((t^{p^m}))$' is true for all $h\in F\setminus f$, 
iff $f$ is a perfect field. All the series $h\in F\setminus f((t^p))$ have a 
polynomial normal form, i.\ e.\ $G_n(h)\bigcap f[t^{-1}, t]\neq \emptyset$. 
All the series $h\in F$ have a polynomial normal form, iff $f$ is a perfect 
field. Every $G_n$-orbit in $F$ is closed. The canonical projection 
\begin{equation} 
\beta: F\longrightarrow F/G_n,\quad h\mapsto G_n(h)
\end{equation}
is closed iff either $f$ is a finite field, or $p=0$. The orbit space 
$F/G_n$ can be metrized by a complete metrics.
\end{theorem}

\begin{proof}
Start with a counterexample to the existence of a polynomial normal form 
for the case of imperfect field.  
Let $f$ be an imperfect field, $c\in f\setminus f^p$ and 
\begin{equation}\label{con}
h=\sum_{i=1}^{\infty} c^{p_i}t^{p^i}
\end{equation}
where
\begin{equation}
p_i=1+p+p^2+\dots+p^i=\frac{p^{i-1}-1}{p-1}.
\end{equation}
Since $h$ satisfies 
\begin{equation}
h-ch^p=t^p ,
\end{equation}
one has 
\begin{equation}\label{coe}
g(h)-cg(h)^p=g(t)^p
\end{equation}     
for any $g\in G_n$. Suppose that $g(h)$ is a polynomial of degree $d$. 
Then $d\geq p$, because as we noticed, $g$ saves the valuation, and 
$\ord h=p$. Calculating coefficients at $t^{d^p}$ in \eqref{coe}, we obtain 
\begin{equation}
c=-\left(\frac{g(t)_d}{g(h)_d}\right)^p \in f^p , 
\end{equation}
a contradiction. Thus, the series \eqref{con} doesn't have a polynomial normal 
form for any imperfect field $f$. 

The next interesting fact, that all the orbits are closed, follows from a 
minimality principle, one of the formulations of Hilbert's 
theorem about bases of polynomial rings, and the following  

\begin{lemma}\label{lem2}
Let $K$ be an algebraically closed field containing $f$; $\A^m$ an affine $K$-space 
of a finite dimension $m$, and $\A^m(f)$ the set of $f$-rational points of $\A^m$. 
For each polynomial function $\Phi: \A^m\rightarrow f$ and for each $f$-linear 
subspace $L$ of the $f$-linear space $K$, one can find an $f$-closed affine 
algebraic variety $S\subseteq \A^m$ such that $S(f)=\Phi^{-1}(L)\bigcap\A^m(f)$.
\end{lemma}

In calculations involving the field $F$, the following lemma is extremely useful.

\begin{lemma}\label{lem3}
Let $p>0\thickspace$, $k=^p\overline{\dots k_i\dots k_1k_0}\in \ZZ_p$, a $p$-adic integer, 
$q_1=^p\overline{\dots q_{1i}\dots q_{11}q_{10}},\thickspace q_2=\\^p\overline
{\dots q_{2i}\dots q_{21}q_{20}}, \dots \in \ZZ_{>0}$, a finite sequence of 
nonnegative integers, $q=q_1+q_2+\dots$. Then
\begin{gather}\label{g1}
\binom{k}{q_1,q_2,\dots}\in \ZZ_p ,\\
\binom{k}{q_1,q_2,\dots}\equiv\prod_{i=0}^{\infty}\binom{k_i}{q_{1i},q_{2i},\dots}
\bmod p ,\label{g2}\\
\binom{k}{q_1,q_2,\dots}\not\equiv 0\bmod p\quad \text{iff 
$\quad \forall i,\thickspace k_i\geq 
q_{1i}+q_{2i}+\dots$ .}\label{last}
\end{gather}
For \eqref{last}, one has $\forall i,\thinspace \nu_p(q_i)\geq \nu_p(k)$ where 
$\nu_p$ denotes the $p$-adic valuation.  
\end{lemma} 

\begin{proof} If $k\in \ZZ_{>0}$, a positive integer, then \eqref{g1} is clear, 
\eqref{g2} follows from the particular case of binomial coefficients, which 
is well-known, and \eqref{last} follows from \eqref{g2}. Since $\ZZ$ is dense 
in $\ZZ_p$, we can extend 
\eqref{g1}, \eqref{g2} and \eqref{last} to $k\in\ZZ_p$ by continuity. 
\end{proof}

By the way, we obtained all the formulas \eqref{g1}, \eqref{g2} and \eqref{last} for 
negative integers $k$ as well, just from the case of positive $k$ and density of 
$\ZZ$ in $\ZZ_p$. For me, that is a very interesting $p$-adic trick. 

The last sentence of Theorem \ref{thm9} follows from  

\begin{lemma}\label{lem4}
Suppose that a group $\Gamma$ acts on a metric space $(X,d)$ by isometries such that 
all the orbits are closed. Then \newcounter{llem4}\begin{list}{(\roman{llem4})} 
{\usecounter{llem4}}
\item\label{i} Function $D:\thickspace (X/\Gamma)^2\rightarrow\RR,\thickspace 
(A,B)\mapsto \inf_{(a,b)\in A\times B} d(a,b)$ is a metric on $X/\Gamma$ defining 
the quotient topology. 
\item\label{ii} $D(A,B)=d(a,B)$ for any $a\in A,\thickspace (A,B)\in (X/\Gamma)^2$.
\item\label{iii} A canonical injection $\inj: X/\Gamma\rightarrow H(X)$ where 
$(H(X),d_H)$ is the space of non-empty closed subsets of $X$ with Hausdorf's 
metrics 
\begin{equation}
d_H(A,B)=\max \{ \sup_{a\in A}d(a,B), \sup_{b\in B} d(A,b)\} ,
\end{equation}
is an isometry and $\inj(X/\Gamma)$ is a closed subset of $H(X)$.
\item\label{iv} If $(X,d)$ is a complete metric space, then $(X/\Gamma, D)$ is 
a complete metric space as well.
\end{list}      
\end{lemma} 
\end{proof}

Theorem \ref{thm9} shows that the situation for $p>0$ is much more complicated 
than it is for $p=0$ which was studied in section \ref{sec2}. In addition to 
invariants described previously, there are a lot of new invariants for $p>0$. For 
instance, the following functions $F\rightarrow \ZZ\bigcup\{\infty\}$ are invariant:
\begin{gather}
\ord_0(h)=\ord(h-x_0(h)) ,\\
\md(h)=\max\{ m\in\ZZ | h\in f((t^{p^m}))\} ,\\
\ord_{\md}(h)=\max \ord\{ h-a | a\in f((t^{p^{\md(h)+1}}))\} ,
\end{gather}
where
\begin{equation}
h=\sum_{i=\ord h}^{\infty} x_i(h) t^i
\end{equation}
with $x_i(h)\in f$ and $x_{\ord h}(h)\neq 0$. 

The next example is more interesting. Determine $w: F\rightarrow\ZZ$, 
\begin{equation}\label{22}
w(h)=\max\left\{ \left[ \frac{\ord_{\md}(h)-m}{|m|_p^{-1}-|\ord_{\md}(h)|_p^{-1}}
\right]\thickspace\left|\right.\thickspace x_m(h)\neq 0, 0\neq m<\ord_{\md}(h)
\right\}
\end{equation}
for $\ord_{\md}(h)>\ord_0(h)$, or $w(h)=0$ otherwise. Here 
\begin{equation}
|m|_p^{-1}\stackrel{\rm def}{=}\max\left\{ p^d\thickspace | \thickspace 
\frac{m}{p^d}\in\ZZ, d\in\ZZ
\right\}
\end{equation}
Letter $w$ in \eqref{22} is the first letter of the word {\em width}.

\begin{prop}
$w$ is $G_0$-invariant.
\end{prop}

\section{An application to QFT}\label{sec5}

Polynomials $P_{mk}^{(-1/\lambda)}$ describing the orbits and invariants of formal 
tensor fields on a line, naturally appear in some other fields of mathematics as 
well. Here is just an example. 

According to \cite{Yu}, denote 
\begin{gather}
P_2=u_2-u_1^2,\\
P_{k+1}=\frac{1}{k+2}\left(\sum_{i=1}^k((i+2)u_{i+1}-2u_1u_i)\frac{\partial P_k}
{\partial u_i}-2ku_1P_k-\sum_{i=2}^{k-1}P_iP_{k+1-i}\right)
\end{gather}
for $k>2$. 

\begin{theorem}\label{thm10}
For $k>1$, 
\begin{equation}
P_k(u_1,u_2,\dots,u_k)=\frac{1}{1-k}P_{0k}^{(1-k)}(1,u_1,u_2,\dots,u_k) .
\end{equation}
\end{theorem}

\begin{proof}
As usual in analogous cases, after guessing the answer, the proof can be done 
by induction on $k$. Limited by the space and time, I omit superfluous details.
\end{proof}

\begin{acn}
I would like to thank Alexandre Kirillov and Fan Chung Graham, as well as my mother  
and my Beautiful and Wonderful wife, Bette.
\end{acn}

\end{document}